\documentclass[a4paper,12pt]{article}

\usepackage{graphicx}
\usepackage{amsmath,amsthm,amssymb,amscd}
\newtheorem{theo}{Theorem}[section]

\newtheorem{lem}{Lemma}[section]
\newtheorem{rem}{Remark}[section]

\newtheorem{defi}{Definition}[section]
\newcommand{\be}{\begin{equation}}
\newcommand{\ee}{\end{equation}}
\newcommand\bes{\begin{eqnarray}} \newcommand\ees{\end{eqnarray}}
\newcommand{\bess}{\begin{eqnarray*}}
\newcommand{\eess}{\end{eqnarray*}}

\begin{document}
\setlength{\baselineskip}{16pt} \pagestyle{myheadings}

\title{Traveling wave solutions for delayed reaction-diffusion systems
\thanks{The work is partially supported by PRC grant NSFC
10671172 and also by the NSF of Jiangsu Province (BK2006064).}
\date{\empty}
\author{  Canrong Tian$^{a}$,  Zhigui Lin$^{b}$ \\
{\small $^a$Basic Department Yancheng Institute of Technology ,
Yancheng 224003, China}\\
{\small $^b$School of Mathematical Science, Yangzhou University,Yangzhou 225002, China}\\
{\small Email: unfoxeses@yahoo.com.cn}}} \maketitle
\begin{quote}
\noindent {\bf Abstract.} {\small This paper is concerned with the
traveling waves of  delayed reaction-diffusion systems where the
reaction function possesses the mixed quasimonotonicity property. By
the so-called monotone iteration scheme and Schauder's fixed point
theorem, it is shown that if the system has a pair of coupled upper
and lower solutions, then there exists at least a traveling wave
solution. More precisely, we reduce the existence of traveling waves
to the existence of an admissible pair of coupled quasi-upper and
quasi-lower solutions which are easy to construct in practice.}

  \vspace{0.5cm}

\end{quote}
\newcommand\HI{{\bf I}}

\section{ Introduction}
Reaction diffusion system is used to model the spatial-temporal pattern. In the past decades,
the traveling wave solutions of the reaction diffusion systems, which are studied as a paradigm for behavior,
have been widely investigated due to significant applications in chemical engineering, population
dynamics and biological models. Since the first instances in which  traveling wave solutions were investigated
were given in 1937 by Kolmogorov et al. \cite{kol} and Fisher \cite{fis}, many methods have been used to study
the traveling wave solutions of various parabolic equations and systems, for example,
 the phase plane technique in \cite{dun1, dun2, kan, kap, mur, tan, tro, vuu},
the degree theory method and the conley index method developed
in \cite{con, gar1, gar2, vol}.

In many realistic models, the delays should be incorporated into the
reaction diffusion system. Due to the presence of delays in the
reaction diffusion system, the classical phase plane technique can
not generate a monotone flow which ensures the existence of the
traveling wave solution. Recently the classical monotone iteration
technique was first used by Wu and Zou \cite{wu, WZ} to establish
the existence of traveling wave solution for delayed reaction
diffusion system. They employed the idea of upper and lower
solutions and an iteration scheme to construct a monotone sequence
of upper solutions which was proved to converge to a solution of the
corresponding wave equation of the reaction-diffusion system under
consideration (see also \cite{huang1,huang2}). In fact, many
researchers had used the monotone iteration technique to prove the
existence of the reaction diffusion system in
\cite{ahm,ker,kun,lad,lak,lee,mar1,mar2,pao1}. In \cite{ma,ma2}, Ma
et al. proved some existence results for traveling wavefronts of
reaction-diffusion systems by using Schauder's fixed point theorem.
One important feature of Ma's method, which was different from the
work of Zou and Wu \cite{wu}, was that the upper solution of the
wave equation was not necessary to converge to two distinct trivial
solutions when $t\rightarrow -\infty$ and $t\rightarrow +\infty$
respectively.  Li et al. developed a new cross iteration scheme and
established the existence of the traveling wave for Lotka-Volterra
competition system with delays \cite{li,li2,li3}.

More recently, Boumenir and Nguyen discussed  in \cite{bou} a
modified version of Perron Theorem for $C^1$-solutions, and set up a
rigorous framework for the monotone iteration method and then apply
it to the predator-prey and Belousov-Zhabotinskii models with
delays.

However, in the iteration process by the monotone iteration it is
required that the nonlinear reaction function possesses a
quasimonotone property in the sector between the upper and lower
solutions \cite{bou, ma, wu}. This paper focus on the delayed
reaction diffusion system without quasimonotonicity.  Motivated by
the above work and the upper and lower solution method developed by
Pao \cite{pao1,pao2,pao3,pao4} and Li et al. \cite{li} for reaction
diffusion systems, we use the coupled upper and lower solutions to
deal with the non-quasimonotonicity, which was first given out in
\cite{lin}. Via the coupled upper and lower solutions, we construct
an appropriate closed bounded convex set. By use of the Schauder 's
fixed point theorem in the convex set, we show the existence of the
traveling wave solution.  Moreover we reduce the existence of
traveling wave solution to the existence of an admissible pair of
quasi-upper and quasi-lower solutions which are easy to construct in
practice.

This paper is organized as follows. In Section 2, we show the
existence of the traveling wave solution by constructing the
classical coupled upper and lower solutions. In Section 3  we relax
the classical coupled upper and lower solutions to the $C^1$ smooth
coupled quasi-upper and quasi-lower solutions. Section 4 deals with
systems with quasimonotone nondecreasing functions, and the
definition of ordered quasi-upper and quasi-lower solutions is
introduced and an existence result of a traveling wavefront is given
by the monotone iteration method. In Section 5 the main result is
illustrated by and applied to a delayed Belousov-Zhabotinskii
equation and a Mutualistic Lotka-Volterra model. This paper ends
with a short discussion.

\section{Coupled upper and lower solutions}
\setcounter{equation}{0}

In this paper, we will consider the following system of
reaction-diffusion systems with time delays \bes
\frac{\partial}{\partial t}\mathbf{u}(x,t)-D
\frac{\partial^2}{\partial^2
x}\mathbf{u}(x,t)=\mathbf{f}(\mathbf{u},
\mathbf{u}_{\mathbf{\tau}}),\label{b1} \ees where $x\in \mathbb{R},$
$t\in (0,\infty),$ $\mathbf{u}\equiv (u_1, \cdots,
u_n)\in\mathbb{R}^n$, $\mathbf{u}_{\mathbf{\tau}}\equiv (u_1(x,
t-\tau_1), \cdots, u_n(x, t-\tau_n))\in\mathbb{R}^n$ for some
positive constants $\tau_1, \cdots, \tau_n$, which are so-called
discrete delays.and $D=diag(d_1,\cdots,d_n)$ with $d_i>0$,
$\mathbf{f}: ~\mathbb{R}^n\times
\mathbb{R}^n\rightarrow\mathbb{R}^n$ is Lipschitz continuous.

For convenience, we denote by $C_b(\mathbb{R}, \mathbb{R})$ the space of all bounded and continuous functions $h : \mathbb{R}\to \mathbb{R}$ endowed
with the super-norm. Moreover, for any $k\in \mathbb{R}_+$, we denote by $C_b^k(\mathbb{R}, \mathbb{R})$
 the space of all continuous differentiable up to the $[k]$-order functions $h$
 such that $d^\gamma h\in C_b(\mathbb{R}, \mathbb{R})$ for any
$|\gamma|\leq [k]$ ($[k]$ denoting the integer part of $k$).
The above spaces for vector-valued functions (with $n$ components) are denoted by $\mathcal{C}_b(\mathbb{R}, \mathbb{R}^n)$ and $\mathcal{C}_b^k(\mathbb{R}, \mathbb{R}^n)$, respectively.

A traveling wave solution of (\ref{b1}) is a special translation
invariant solution of the form
$\mathbf{u}(x,t)=\mathbf{\varphi}(t+x/c)$, where
$\mathbf{\varphi}\in \mathcal{C}^2_b(\mathbb{R},\mathbb{R}^n)$ is
the profile of the wave and $c>0$ is a constant corresponding to the
wave speed. The vibration at the space point $x=0$ is
$\mathbf{u}(t)=\mathbf{\varphi}(t)$; The vibration $\mathbf{u}(t)$
propagating from the space value $x=0$ to $x$ costs the time $x/c$,
where $c$ is wave velocity, in the case $c>0$ the traveling wave
move to the left, in the case $c< 0$ the traveling wave move to the
right. Our definition is more visual than
$\mathbf{u}(x,t)=\mathbf{\varphi}(x+ct)$ in \cite{bou,ma,wu}.

 Substituting $\mathbf{u}(x,t)=\mathbf{\varphi}(t+x/c)$ into
(\ref{b1}) and letting $s=t+x/c$, denoting also $t$, we obtain
the\bes \mathbf{\varphi}'(t)-\frac{D}{c^2}
\mathbf{\varphi}''(t)=\mathbf{f}(\mathbf{\varphi}(t),
\mathbf{\varphi}(t-\mathbf{\tau})),\quad t\in\mathbb{R}. \label{b2}
\ees corresponding wave equations  If for some wave velocity $c$,
(\ref{b2}) has a solution $\mathbf{\varphi}$ defined on $\mathbb{R}$
such that \bes \lim_{t\rightarrow
-\infty}\mathbf{\varphi}(t)=\mathbf{u}_{-},\quad \lim_{t\rightarrow
+\infty}\mathbf{\varphi}(t)=\mathbf{u}_{+} \label{b3} \ees exist,
then $\mathbf{u}(x,t)=\mathbf{\varphi}(t+x/ct)$ is called traveling
wave with speed $c$. Moreover, if $\mathbf{\varphi}$ is monotone in
$t\in \mathbb{R}$, then it is called a traveling wavefront.

Without loss of generality, we can assume
$\mathbf{u}_{-}=0$ and $\mathbf{u}_{+} =\mathbf{K}>0$. Let \bess
\mathcal{C}_{[0,\mathbf{K}]}(\mathbb{R},\mathbb{R}^n)=\{\mathbf{\varphi}=(\varphi_1(t), \cdots, \varphi_n(t))\in
\mathcal{C}_b(\mathbb{R}, \mathbb{R}^n):\,  0\leq {\varphi}_i(t)\leq\ K_i,\, t\in\mathbb{R}\}.
\eess Our aim is looking for a solution of (\ref{b2}) in
$\mathcal{C}_{[0,\mathbf{K}]}(\mathbb{R}, \mathbb{R}^n)$. Throughout this
paper, the following hypothesis will be imposed on the reaction term $\mathbf{f}$:
\begin{itemize}
  \item [$(H_1)$]$\mathbf{f}(0)=\mathbf{f}(\mathbf{K})=0.$
\end{itemize}
Obviously, we should replace (\ref{b3}) with \bes \lim_{t\rightarrow
-\infty}\mathbf{\varphi}(t)=0,\quad \lim_{t\rightarrow
+\infty}\mathbf{\varphi}(t)=\mathbf{K}\label{b4} \ees In this paper,
we explore the existence of the traveling wave solutions of
(\ref{b1}) where the reaction term $\mathbf{f}$ is mixed
quasimonotone.
\begin{itemize}
  \item [$(H_2)$]The function $\mathbf{f}(\mathbf{u}, \mathbf{u}_{\mathbf{\tau}})=(f_1(\mathbf{u}, \mathbf{u}_{\mathbf{\tau}}),\cdots, f_n(\mathbf{u}, \mathbf{u}_{\mathbf{\tau}}))$
   is  a $C^1$ function and possesses a mixed quasimonotone property in a
  subset $[0,\mathbf{K}]$ of $\mathbb{R}^n$. \end{itemize}
  The above hypothesis implies that there exist constants $\beta_i$ such that $f_i$ satisfies the Lipschity condition
  $$|f_i(\mathbf{u}, \mathbf{u}_{\mathbf{\tau}})- f_i(\mathbf{v}, \mathbf{v}_{\mathbf{\tau}})|\leq \beta_i(||\mathbf{u}-\mathbf{v}||+
  ||\mathbf{u}_{\mathbf{\tau}}-\mathbf{v}_{\mathbf{\tau}}||)$$
  for all $\mathbf{u}, \mathbf{v}, \mathbf{u}_{\mathbf{\tau}}$ and $\mathbf{v}_{\mathbf{\tau}}$ in
  $\mathcal{C}_{[0,\mathbf{K}]}(\mathbb{R},\mathbb{R}^n)$, $i=1, \cdots, n$, where $|\cdot|$ and $||\cdot||$ denote the super norm in $\mathbb{R}^n$ and $\mathcal{C}(\mathbb{R}, \mathbb{R}^n)$, respectively.

  Recall that  by writing vectors $\mathbf{u}$ and $\mathbf{u}_{\mathbf{\tau}}$ in $\mathbb{R}^n$ in the split form
  \bess
\mathbf{u}\equiv(u_i,[\mathbf{u}]_{a_i},[\mathbf{u}]_{b_i}),\ \mathbf{u}_{\mathbf{\tau}}\equiv([\mathbf{u}_{\mathbf{\tau}}]_{c_i},[\mathbf{u}_{\mathbf{\tau}}]_{d_i})
  \eess
  where $[\mathbf{u}]_{\sigma}$ denotes a vector with $\sigma$
  components of $\mathbf{u}$, the function $\mathbf{f}(\mathbf{u}, \mathbf{u}_{\mathbf{\tau}})$
  is said to have a mixed quasimonotone property if for each
  $i=1,\cdots,n,$ there exist nonnegative integers $a_i$, $b_i$, $c_i$ and $d_i$ with
  $a_i+b_i=n-1$ and $c_i+d_i=n$ such that $f_i(\mathbf{u}, \mathbf{v})$ is monotone
  nondecreasing in $[\mathbf{u}]_{a_i}$ and $[\mathbf{v}]_{c_i}$ and is monotone
  nonincreasing in $[\mathbf{u}]_{b_i}$ and $[\mathbf{v}]_{d_i}$. If $b_i=d_i=0$ for all $i$ then
  $\mathbf{f}(\mathbf{u}, \mathbf{v})$ is said to be quasimonotone nondecreasing.

  The above general assumptions are used to establish the existence of
  a traveling wave solution to (\ref{b2}). Our approach to the
  problem is by the method of coupled upper and lower solutions which are
  defined as follows:
  \begin{defi}
A pair of vectors $\mathbf{\tilde\varphi}\equiv(\tilde
\varphi_1,\cdots,\tilde \varphi_n),$
$\mathbf{\hat\varphi}\equiv(\hat \varphi_1,\cdots,\hat \varphi_n)$
in $\mathcal{C}^2_b(\mathbb{R}, \mathbb{R}^n)$ are called coupled upper and lower
solutions of (\ref{b2}) if
$\mathbf{\tilde\varphi}\geq\mathbf{\hat\varphi}$ and if \bes
\begin{array}{ll}
\tilde\varphi'_i(t)-\frac{d_i}{c^2}\tilde\varphi''_i(t)\geq
f_i(\tilde\varphi_i,~[\mathbf{\tilde\varphi}]_{a_i},~[\mathbf{\hat\varphi}]_{b_i},
~[\mathbf{\tilde\varphi}_{\mathbf{\tau}}]_{c_i},~[\mathbf{\hat\varphi}_{\mathbf{\tau}}]_{d_i}),\\
\hat\varphi'_i(t)-\frac{d_i}{c^2}\hat\varphi''_i(t)\leq
f_i(\hat\varphi_i,~[\mathbf{\hat\varphi}]_{a_i},~[\mathbf{\tilde\varphi}]_{b_i},
~[\mathbf{\hat\varphi}_{\mathbf{\tau}}]_{c_i},~[\mathbf{\tilde\varphi}_{\mathbf{\tau}}]_{d_i})&(i=1,\cdots,n),
\end{array}\label{b5}
\ees where
$\mathbf{\varphi}_{\mathbf{\tau}}(t)=\mathbf{\varphi}(t-\mathbf{\tau})$.
\end{defi}
\begin{rem}
If $\mathbf{f}(\mathbf{\varphi}, \mathbf{\varphi}_{\mathbf{\tau}})$ is quasimonotone nondecreasing, that is,
$b_i=d_i=0$ for all $i$,  then the pair of vectors called ordered upper and lower
solutions of (\ref{b2})
\end{rem}

Since that $\mathbf{f}$ satisfies $(H_2)$ and Lipschitz continuous,
we have \bes\begin{split}& f_i(\varphi_i,
~[\mathbf{\varphi}]_{a_i},~[\mathbf{\varphi}]_{b_i},
~[\mathbf{\varphi}_{\mathbf{\tau}}]_{c_i},~[\mathbf{\varphi}_{\mathbf{\tau}}]_{d_i})-
f_i(\varphi'_i, ~[\mathbf{\varphi}]_{a_i},~[\mathbf{\varphi}]_{b_i},
~[\mathbf{\varphi}_{\mathbf{\tau}}]_{c_i},~[\mathbf{\varphi}_{\mathbf{\tau}}]_{d_i})\\&+
\beta_i(\varphi_i-\varphi'_i)\geq 0 \ \textrm{ for all }\  0\leq
\varphi'_i\leq\varphi_i\leq K_i, \quad
i=1,\cdots,n.\end{split}\label{b6}\ees Next we define an operator
$\mathbf{H}:~\mathcal{C}_{[0,\mathbf{K}]}(\mathbb{R},\mathbb{R}^n)\rightarrow
\mathcal{C}_{[0,\mathbf{K}]}(\mathbb{R},\mathbb{R}^n)$ by \bes
\mathbf{H}(\mathbf{\varphi}, \mathbf{\varphi}_{\mathbf{\tau}})(t)
=\mathbf{f}(\mathbf{\varphi},
\mathbf{\varphi}_{\mathbf{\tau}})+\beta\mathbf{\varphi}(t), \quad
\mathbf{\varphi}\in
C_{[0,\mathbf{K}]}(\mathbb{R},\mathbb{R}^n),\label{b7}\ees where
$\mathbf{H}=(H_1,\cdots,H_n),$ $H_i(\mathbf{\varphi},
\mathbf{\varphi}_{\mathbf{\tau}})(t)=f_i(\mathbf{\varphi},
\mathbf{\varphi}_{\mathbf{\tau}})+\beta_i\varphi_i(t).$ Clearly,
with the above notations, the system (\ref{b2}) is equivalent to the
following system of ordinary differential equations \bes
\mathbf{\varphi}'(t)-\frac{D}{c^2}
\mathbf{\varphi}''(t)+\beta\mathbf{\varphi}(t)=\mathbf{H}(\mathbf{\varphi},
\mathbf{\varphi}_{\mathbf{\tau}})(t),\quad t\in\mathbb{R}.
\label{b8}\ees Our first iteration involves the following linear
system of ordinary differential equations \bes
\begin{array}{ll}
c(\overline x_i^{(1)})'-\frac{d_i}{c^2}(\overline
x_i^{(1)})''+\beta_i\overline x_i^{(1)}=\beta_i\tilde \varphi_i+
f_i(\tilde\varphi_i,~[\mathbf{\tilde\varphi}]_{a_i},~[\mathbf{\hat\varphi}]_{b_i},
~[\mathbf{\tilde\varphi}_{\tau}]_{c_i},~[\mathbf{\hat\varphi}_{\tau}]_{d_i}),\\
c(\underline x_i^{(1)})'-\frac{d_i}{c^2}(\underline
x_i^{(1)})''+\beta_i\underline x_i^{(1)}=\beta_i\hat \varphi _i+
f_i(\hat\varphi_i,~[\mathbf{\hat\varphi}]_{a_i},~[\mathbf{\tilde\varphi}]_{b_i},
~[\mathbf{\hat\varphi}_{\tau}]_{c_i},~[\mathbf{\tilde\varphi}_{\tau}]_{d_i}).&
\end{array}\label{b9}
\ees Note that \bess \lambda_{1i}=\frac{c^2(1-\sqrt{1+4\beta_i
d_i/c^2})}{2d_i},\quad \lambda_{2i}=\frac{c^2(1+\sqrt{1+4\beta_i
d_i/c^2})}{2d_i}\eess are the negative and positive real roots of
the equation
$$\frac{d_i}{c^2}\lambda^2-\lambda-\beta_i=0,\ \ i=1, 2, \cdots, n.$$
Using the Perron Theorem yields
 \bes\begin{split} \overline
x_i^{(1)}=\frac{c^2}{d_i(\lambda_{2i}-\lambda_{1i})}(\int_{-\infty}^t
e^{\lambda_{1i}(t-s)}(\beta_i\tilde\varphi_i+f_i(\tilde\varphi_i,
[\mathbf{\tilde\varphi}]_{a_i}, [\mathbf{\hat\varphi}]_{b_i},
[\mathbf{\tilde\varphi}_{\tau}]_{c_i},
[\mathbf{\hat\varphi}_{\tau}]_{d_i}))ds\\+\int_{t}^{+\infty}
e^{\lambda_{2i}(t-s)}(\beta_i\tilde\varphi_i+f_i(\tilde\varphi_i,
[\mathbf{\tilde\varphi}]_{a_i},
[\mathbf{\hat\varphi}]_{b_i}, [\mathbf{\tilde\varphi}_{\tau}]_{c_i}, [\mathbf{\hat\varphi}_{\tau}]_{d_i}))ds),\\
\underline
x_i^{(1)}=\frac{c^2}{d_i(\lambda_{2i}-\lambda_{1i})}(\int_{-\infty}^t
e^{\lambda_{1i}(t-s)}(\beta_i\hat\varphi_i+f_i(\hat\varphi_i,
[\mathbf{\hat\varphi}]_{a_i}, [\mathbf{\tilde\varphi}]_{b_i},
[\mathbf{\hat \varphi}_{\tau}]_{c_i},
[\mathbf{\tilde\varphi}_{\tau}]_{d_i}))ds\\+\int_{t}^{+\infty}
e^{\lambda_{2i}(t-s)}(\beta_i\hat\varphi_i+f_i(\hat\varphi_i,
[\mathbf{\hat\varphi}]_{a_i}, [\mathbf{\tilde\varphi}]_{b_i},
[\mathbf{\hat\varphi}_{\tau}]_{c_i},
[\mathbf{\tilde\varphi}_{\tau}]_{d_i}))ds)
\end{split}\label{b10}\ees
for $i=1, \cdots, n$. Then by Lemma 2.1 of \cite{lin},
$\mathbf{\overline x}^{(1)}\equiv (\overline
x_1^{(1)},\cdots,\overline x_n^{(1)})$ and $\mathbf{\underline
x}^{(1)}\equiv (\underline x_1^{(1)},\cdots,\underline x_n^{(1)})$
have the following properties.
\begin{lem}
Let $\mathbf{\overline x}^{(1)}$ and $\mathbf{\underline x}^{(1)}$
be the solution of (\ref{b9}), then we have
\begin{itemize}
   \item [(i)] $\mathbf{\hat \varphi}\leq\mathbf{\underline x}^{(1)}\leq
  \mathbf{\overline x}^{(1)}\leq\mathbf{\tilde \varphi} $.
  \item [(ii)]$\mathbf{\overline x}^{(1)}$, $\mathbf{\underline
  x}^{(1)}$ are a pair of  coupled upper and lower solutions of
  (\ref{b2}).
\end{itemize}\label{lem11}
\end{lem}

Now in order to prove the existence of the traveling wave solution,
we are in the position to apply the  Schauder's fixed point theorem.
We define the operator
\bess\mathbf{F}=(F_1,\cdots,F_n):\mathcal{C}_{[0,\mathbf{K}]}(\mathbb{R},\mathbb{R}^n)\rightarrow\mathcal{C}_{[0,\mathbf{K}]}(\mathbb{R},\mathbb{R}^n)\eess
by \bess
&(F_i\varphi_i)(t)\\=&\frac{c^2}{d_i(\lambda_{2i}-\lambda_{1i})}(\int_{-\infty}^t
e^{\lambda_{1i}(t-s)}H_i(\mathbf{\varphi},
\mathbf{\varphi}_{\mathbf{\tau}})(s)ds+\int_{t}^{+\infty}
e^{\lambda_{2i}(t-s)}H_i(\mathbf{\varphi},
\mathbf{\varphi}_{\mathbf{\tau}})(s)ds) \eess for $i=1,2,\cdots,n.$

Let $\rho>0$ be such that
$\rho<\min\{-\lambda_{1i},\lambda_{2i}:~i=1,2,\cdots,n\}$, and let
\bess
&\mathbf{B}_{\rho}(\mathbb{R},\mathbb{R}^n)=\{\mathbf{\varphi}\in\mathcal{C}_{[0,\mathbf{K}]}(\mathbb{R},\mathbb{R}^n)
:~sup_{t\in\mathbb{R}}|\mathbf{\varphi}(t)|e^{-\rho|t|}<\infty\},\\&
|\mathbf{\varphi}|_\rho=\sup_{t\in\mathbb{R}}|\mathbf{\varphi}(t)|e^{-\rho|t|}.
\eess Then it is easy to check that
$\mathbf{B}_{\rho}(\mathbb{R},\mathbb{R}^n),|\cdot|_\rho$ is a
Banach space.

\begin{lem}
Let the closed convex set $\mathbf{\Gamma}=\{\mathbf{\varphi}\in
\mathcal{C}_{[0,\mathbf{K}]}(\mathbb{R},\mathbb{R}^n):~\mathbf{\hat\varphi}\leq\mathbf{\varphi}\leq\mathbf{\tilde\varphi}\}$,
where $\mathbf{\tilde\varphi}$ and $\mathbf{\hat\varphi}$ are
coupled upper and lower solutions of (\ref{b2}), then
$\mathbf{F}(\mathbf{\Gamma})\subseteq \mathbf{\Gamma}$.\label{lem1}
\end{lem}
\begin{proof}
Since that the mixed quasimonotone property  of the operator
$\mathbf{H}$, $\forall \mathbf{\varphi}\in\mathbf{\Gamma}$, we have
\bess (F_i\varphi_i)(t)\leq
\frac{c^2}{d_i(\lambda_{2i}-\lambda_{1i})}(\int_{-\infty}^t
e^{\lambda_{1i}(t-s)}(\beta_i\tilde\varphi_i+f_i(\tilde\varphi_i,
[\mathbf{\tilde\varphi}]_{a_i}, [\mathbf{\hat\varphi}]_{b_i},
[\mathbf{\tilde\varphi}_{\tau}]_{c_i},
[\mathbf{\hat\varphi}_{\tau}]_{d_i}))ds\\+\int_{t}^{+\infty}
e^{\lambda_{2i}(t-s)}(\beta_i\tilde\varphi_i+f_i(\tilde\varphi_i,
[\mathbf{\tilde\varphi}]_{a_i}, [\mathbf{\hat\varphi}]_{b_i},
[\mathbf{\tilde\varphi}_{\tau}]_{c_i},
[\mathbf{\hat\varphi}_{\tau}]_{d_i}))ds)\eess By the virtue of Lemma
\ref{lem11}, we induce
$(F_i\varphi_i)(t)\leq\mathbf{\tilde\varphi}_i(t)$. Similarly, we
have $(F_i\varphi_i)(t)\geq\mathbf{\hat\varphi}_i(t)$ for all
$i=1,\cdots,n$. Therefore
$\mathbf{F}(\mathbf{\Gamma})\subseteq\mathbf{\Gamma}$.
\end{proof}

\begin{lem}
If the hypothesis $(H_2)$ holds, then
$\mathbf{F}:\mathcal{C}_{[0,\mathbf{K}]}(\mathbb{R},\mathbb{R}^n)\rightarrow\mathcal{C}_{[0,\mathbf{K}]}(\mathbb{R},\mathbb{R}^n)$
is continuous with respect to the norm $|\cdot|_\rho$ in
$\mathbf{B}_{\rho}(\mathbb{R},\mathbb{R}^n)$.\label{lem2}
\end{lem}
\begin{proof}
$\forall
\mathbf{\varphi},~\mathbf{\varphi'}\in\mathcal{C}_{[0,\mathbf{K}]}(\mathbb{R},\mathbb{R}^n),$
in view of the definition of $\mathbf{F}$,  we have
\bes\begin{split}
F_i(\mathbf{\varphi},\mathbf{\varphi}_\mathbf{\tau})-
F_i(\mathbf{\varphi}',\mathbf{\varphi}'_\mathbf{\tau})=\frac{c^2}{d_i(\lambda_{2i}-\lambda_{1i})}(\int_{-\infty}^te^{\lambda_{1i}(t-s)}(\beta_i(\varphi_i
-\varphi'_i)+f_i(\mathbf{\varphi},\mathbf{\varphi}_\mathbf{\tau})\\-f_i(\mathbf{\varphi}',\mathbf{\varphi}'_\mathbf{\tau}))ds
+\int_t^{+\infty}e^{\lambda_{2i}(t-s)}(\beta_i(\varphi_i
-\varphi'_i)+f_i(\mathbf{\varphi},\mathbf{\varphi}_\mathbf{\tau})-f_i(\mathbf{\varphi}',\mathbf{\varphi}'_\mathbf{\tau}))ds).\end{split}
\label{bb1}\ees Since that $\mathbf{\varphi},~\mathbf{\varphi'}\leq
\mathbf{K}$, $(H_2)$ implies that
$f_i(\mathbf{\varphi},\mathbf{\varphi}_\mathbf{\tau})$ is bounded
for $\mathbf{\varphi}\in
\mathcal{C}_{[0,\mathbf{K}]}(\mathbb{R},\mathbb{R}^n)$. Then the
term $\beta_i(\varphi_i
-\varphi'_i)+f_i(\mathbf{\varphi},\mathbf{\varphi}_\mathbf{\tau})-f_i(\mathbf{\varphi}',\mathbf{\varphi}'_\mathbf{\tau})$
is bounded, for convenience, we denote $\beta_i(\varphi_i
-\varphi'_i)+f_i(\mathbf{\varphi},\mathbf{\varphi}_\mathbf{\tau})-f_i(\mathbf{\varphi}',\mathbf{\varphi}'_\mathbf{\tau})\leq
M_i$. (\ref{bb1}) is transformed into \bes\begin{split}
F_i(\mathbf{\varphi},\mathbf{\varphi}_\mathbf{\tau})-
F_i(\mathbf{\varphi}',\mathbf{\varphi}'_\mathbf{\tau})\leq&\frac{c^2M_i}{d_i(\lambda_{2i}-\lambda_{1i})}(\int_{-\infty}^te^{\lambda_{1i}(t-s)}ds
+\int_t^{+\infty}e^{\lambda_{2i}(t-s)}ds)\\
=&\frac{c^2M_i}{d_i(\lambda_{2i}-\lambda_{1i})}(\frac{1}{\lambda_{2i}}-\frac{1}{\lambda_{1i}})=-\frac{c^2M_i}{d_i\lambda_{1i}\lambda_{2i}}.\end{split}
\label{bb2}\ees It follows from (\ref{bb2}) that
\bes\begin{split}|F_i(\mathbf{\varphi},\mathbf{\varphi}_\mathbf{\tau})-
F_i(\mathbf{\varphi}',\mathbf{\varphi}'_\mathbf{\tau})|
e^{-\rho|t|}\leq-\frac{c^2M_i}{d_i\lambda_{1i}\lambda_{2i}}e^{-\rho|t|}\leq-\frac{c^2M_i}{d_i\lambda_{1i}\lambda_{2i}}.\end{split}
\label{bb3}\ees Therefore
$|F_i(\mathbf{\varphi},\mathbf{\varphi}_\mathbf{\tau})-
F_i(\mathbf{\varphi}',\mathbf{\varphi}'_\mathbf{\tau})|_\rho\leq-\frac{c^2M_i}{d_i\lambda_{1i}\lambda_{2i}}$
for all $i=1,\cdots,n$. That is,
$\mathbf{F}:\mathcal{C}_{[0,\mathbf{K}]}(\mathbb{R},\mathbb{R}^n)\rightarrow\mathcal{C}_{[0,\mathbf{K}]}(\mathbb{R},\mathbb{R}^n)$
is continuous.
\end{proof}

\begin{lem}
If the hypothesis $(H_2)$ holds and $\mathbf{\Gamma}$ is defined in
Lemma \ref{lem1}, then
$\mathbf{F}:\mathbf{\Gamma}\rightarrow\mathbf{\Gamma}$ is
compact.\label{lem3}
\end{lem}
\begin{proof}
First we compute $|\frac{dF_i}{dt}|_\rho$, for any
$\mathbf{\varphi}\in\mathcal{C}_{[0,\mathbf{K}]}(\mathbb{R},\mathbb{R}^n)$,
we have \bes
\begin{split}
\frac{dF_i}{dt}(\mathbf{\varphi},\mathbf{\varphi}_\mathbf{\tau})(t)=\frac{c^2\lambda_{1i}}{d_i(\lambda_{2i}-\lambda_{1i})}
\int_{-\infty}^te^{\lambda_{1i}(t-s)}(\beta_i\varphi_i
+f_i(\mathbf{\varphi},\mathbf{\varphi}_\mathbf{\tau}))ds\\+\frac{c^2\lambda_{2i}}{d_i(\lambda_{2i}-\lambda_{1i})}
\int_t^{+\infty}e^{\lambda_{2i}(t-s)}(\beta_i\varphi_i
+f_i(\mathbf{\varphi},\mathbf{\varphi}_\mathbf{\tau}))ds.
\end{split}\label{bb4}
\ees It follows  from the similar argument of (\ref{bb2}) that
 \bes
\begin{split}
\frac{dF_i}{dt}(\mathbf{\varphi},\mathbf{\varphi}_\mathbf{\tau})(t)\leq\frac{c^2\lambda_{1i}M_i}{d_i(\lambda_{2i}-\lambda_{1i})}
\frac{1}{\lambda_{1i}}+\frac{c^2\lambda_{2i}M_i}{d_i(\lambda_{2i}-\lambda_{1i})}
\frac{1}{\lambda_{2i}}=\frac{2c^2M_i}{d_i(\lambda_{2i}-\lambda_{1i})}.
\end{split}\label{bb5}\ees
It follows from (\ref{bb5}) that
\bes\begin{split}|\frac{dF_i}{dt}(\mathbf{\varphi},\mathbf{\varphi}_\mathbf{\tau})(t)|
e^{-\rho|t|}\leq\frac{2c^2M_i}{d_i(\lambda_{2i}-\lambda_{1i})}e^{-\rho|t|}\leq\frac{2c^2M_i}{d_i(\lambda_{2i}-\lambda_{1i})}.\end{split}
\label{bb6}\ees Therefore
$|\frac{dF_i}{dt}(\mathbf{\varphi},\mathbf{\varphi}_\mathbf{\tau})(t)|_\rho\leq\frac{2c^2M_i}{d_i(\lambda_{2i}-\lambda_{1i})}$
for all $i=1,\cdots,n$. Hence $\mathbf{F}$ is equicontinuous on
$\mathcal{C}_{[0,\mathbf{K}]}(\mathbb{R},\mathbb{R}^n)$. In view of
Lemma \ref{lem1}, $\mathbf{F}(\mathbf{\Gamma})$ is uniformly
bounded.

Next we claim that
$\mathbf{F}:\mathbf{\Gamma}\rightarrow\mathbf{\Gamma}$ is compact.
Define the operator sequence $\{\mathbf{F}^{(n)}\}$, where
$\mathbf{F}^{(n)}:\mathcal{C}_{[0,\mathbf{K}]}(\mathbb{R},\mathbb{R}^n)\rightarrow\mathcal{C}_{[0,\mathbf{K}]}(\mathbb{R},\mathbb{R}^n)$
by \bess
\mathbf{F}^{(n)}(\mathbf{\varphi},\mathbf{\varphi}_\mathbf{\tau})(t)=\left\{\begin{array}{ll}
 \mathbf{F}(\mathbf{\varphi},\mathbf{\varphi}_\mathbf{\tau})(-n),&t\in(-\infty,-n),\\
  \mathbf{F}(\mathbf{\varphi},\mathbf{\varphi}_\mathbf{\tau})(t),&t\in[-n,n], \\
  \mathbf{F}(\mathbf{\varphi},\mathbf{\varphi}_\mathbf{\tau})(n),&t\in(n,+\infty).
   \end{array}\right.\eess
Hence the sequence $\{\mathbf{F}^{(n)}\}$ are uniformly bounded and
equicontinuous. It follows from Arzela-Ascoli theorem that
$\mathbf{F}^{(n)}$ is compact. Therefore we have \bes
\begin{split}&
|F_i^{(n) }-F_i|_\rho=\sup_{t\in
\mathbb{R}}|F_i^n(\mathbf{\varphi},\mathbf{\varphi}_\mathbf{\tau})(t)-F_i(\mathbf{\varphi},\mathbf{\varphi}_\mathbf{\tau})(t)|
e^{-\rho|t|}\\&
=\sup_{t\in(-\infty,-n)\cup(n,+\infty)}|F_i^n(\mathbf{\varphi},\mathbf{\varphi}_\mathbf{\tau})(t)-F_i(\mathbf{\varphi},\mathbf{\varphi}_\mathbf{\tau})(t)|
e^{-\rho|t|}\\& \leq 2K_ie^{-\rho n}\rightarrow0 \texttt{ as
}t\rightarrow\infty,
\end{split}\label{bb7}
\ees where
$\mathbf{\varphi}\in\mathcal{C}_{[0,\mathbf{K}]}(\mathbb{R},\mathbb{R}^n)$.
By the virtue of proposition 2.1 in \cite{zei}, the sequence
$\{\mathbf{F}^{(n)}\}$ converges to $\mathbf{F}$ in
$\mathbf{\Gamma}$ with respect to the norm $|\cdot|_\rho$. Therefore
$\mathbf{F}:\mathbf{\Gamma}\rightarrow\mathbf{\Gamma}$ is compact.

\end{proof}

\begin{theo}
Assume that ($H_1$) and ($H_2$) hold. Suppose that
$\mathbf{\tilde\varphi}\in
\mathcal{C}_{[0,\mathbf{K}]}(\mathbb{R},\mathbb{R}^n)$,
$\mathbf{\hat\varphi}\in
\mathcal{C}_{[0,\mathbf{K}]}(\mathbb{R},\mathbb{R}^n)$ be a pair of
upper and lower solutions of (\ref{b2}), and \bes
\lim_{t\rightarrow-\infty} \mathbf{\tilde\varphi}(t)=0,\quad
\lim_{t\rightarrow+\infty} \mathbf{\hat\varphi}(t)=
\mathbf{K},\label{b18}\ees
 then (\ref{b2}) and (\ref{b4}) admit a solution.
That is, the problem (\ref{b1}) has a traveling wave solution.
\end{theo}
\begin{proof}
First we define the following profile set such as in Lemma
\ref{lem1}\bess \mathbf{\Gamma}=\{\mathbf{\varphi}\in
\mathcal{C}_{[0,\mathbf{K}]}(\mathbb{R},\mathbb{R}^n):~\mathbf{\hat\varphi}\leq\mathbf{\varphi}\leq\mathbf{\tilde\varphi}\},
\eess it is easy to show that $\mathbf{\Gamma}$ is a closed convex
set.

Now we define the operator such as in Lemma \ref{lem2}
\bess\mathbf{F}:~\mathbf{\Gamma}\rightarrow\mathbf{\Gamma},\eess in
view of Lemma \ref{lem1}, Lemma \ref{lem2}, Lemma \ref{lem3},
$\mathbf{F}$ is continuous and compact with respect to
$|\cdot|_\rho$. By Schauder's fixed point theorem, there exists a
fixed point $\mathbf{\varphi^*}\in\mathbf{\Gamma}$ such that
$\mathbf{F}(\mathbf{\varphi^*},\mathbf{\varphi^*_\tau})=\mathbf{\varphi^*}$.
Using the Perron Theorem, (\ref{b8}) has a solution
$\mathbf{\varphi^*}$, that is $\mathbf{\varphi^*}$ is a solution of
(\ref{b2}).

Next we will show that $\mathbf{\varphi^*}$ satisfies the boundary
condition (\ref{b4}). In view of Lemma \ref{lem1},
$\mathbf{\hat\varphi}\leq\mathbf{\varphi^*}\leq\mathbf{\tilde\varphi}$.
It follows from (\ref{b18}) that \bess
0\leq\lim_{t\rightarrow-\infty}\mathbf{\varphi^*}(t)\leq\lim_{t\rightarrow-\infty}\mathbf{\tilde\varphi}=0,\\
\mathbf{K}\geq\lim_{t\rightarrow+\infty}\mathbf{\varphi^*}(t)\geq\lim_{t\rightarrow+\infty}\mathbf{\hat\varphi}=\mathbf{K}.\eess
Therefore $\mathbf{\varphi^*}$ is a traveling wave solution of the
problem (\ref{b1}). Thus the proof is completed.\label{them1}
\end{proof}

\section{Coupled quasi-upper and quasi-lower solutions}
\setcounter{equation}{0}

In Theorem 2.1, we see that the smooth conditions on the coupled
upper and lower solutions are too strong. In fact, it is very
difficult to seek the $C^2$ smooth coupled upper  and lower
solutions for a special model. We intend to relax the smoothness of
the upper and lower solutions to $C^1$. Thus we should cite the
modified Perron theorem in \cite{bou}.
\begin{defi}
Considering the following scalar ordinary equation \bes
u''(t)+\alpha u'(t)+\beta u(t)+f(t)=0,~t\in \mathbb{R},~u(t)\in
\mathbb{R}\label{cc1},\ees where $\beta<0$, $f$ is a bounded and
continuous function on $\mathbb{R}\setminus\{0\}$ and both $f(0^+)$
and $f(0^-)$ exist. Then, a function $u$ defined on $\mathbb{R}$ is
said to be a generalized solution of (\ref{cc1}) if

(1) $u$ and $u'$ are bounded and continuous on $\mathbb{R}$.

(2) $u''$ exists and is continuous on  $\mathbb{R}\setminus\{0\},$
and both $u''(0^-)$ and $u''(0^+)$ exist.
\end{defi}

\begin{lem}(\cite{bou}) Consider (\ref{cc1}) with $\beta<0$, and
assume that

(1) $f$ is a bounded and continuous function on
$\mathbb{R}\setminus\{0\}$ and both $f(0^+)$ and $f(0^-)$ exist,

(2) (\ref{cc1}) holds in the classical sense for all $t$ except
possibly at $t=0$.

Then (\ref{cc1}) has a unique generalized solution $u$ given by \bes
u(t)=\frac{1}{\lambda_2-\lambda_1}(\int_{-\infty}^te^{\lambda_1(t-s)}f(s)ds+\int_t^{+\infty}e^{\lambda_2(t-s)}f(s)ds),\label{cc2}
\ees where $\lambda_1$ and  $\lambda_2$ are respectively the
negative and positive roots of
$\lambda^2+\alpha\lambda+\beta=0$.\label{lemc1}
\end{lem}

 Now we give the following definition of an admissible
upper and lower solutions.
\begin{defi}
Assume that $\mathbf{\tilde\varphi}$, $\mathbf{\hat\varphi}\in
\mathcal{C}^1_b(\mathbb{R},\mathbb{R}^n)$ ,
$\frac{d^2\mathbf{\tilde\varphi}}{dt^2}(t)$ and
$\frac{d^2\mathbf{\hat\varphi}}{dt^2}(t)$ exist and continuous on
$\mathbb{R}\setminus\{0\}$, and\bes\begin{split}
\sup_{t\rightarrow\mathbb{R}\setminus\{0\}}
|\frac{d^2\mathbf{\hat\varphi}}{dt^2}(t)|<+\infty, \texttt{ and
}\lim_{t\rightarrow
0^-}\frac{d^2\mathbf{\hat\varphi}}{dt^2}(t),~\lim_{t\rightarrow
0^+}\frac{d^2\mathbf{\hat\varphi}}{dt^2}(t)\texttt{ exist,
}\\
\sup_{t\rightarrow\mathbb{R}\setminus\{0\}}
|\frac{d^2\mathbf{\tilde\varphi}}{dt^2}(t)|<+\infty, \texttt{ and
}\lim_{t\rightarrow
0^-}\frac{d^2\mathbf{\tilde\varphi}}{dt^2}(t),~\lim_{t\rightarrow
0^+}\frac{d^2\mathbf{\tilde\varphi}}{dt^2}(t)\texttt{ exist,
}\label{c1}\end{split}\ees $\mathbf{\tilde\varphi}$,
$\mathbf{\hat\varphi}$ satisfy \bes
\begin{array}{ll}
\tilde\varphi'_i(t)-\frac{d_i}{c^2}\tilde\varphi''_i(t)\geq
f_i(\tilde\varphi_i,~[\mathbf{\tilde\varphi}]_{a_i},
~[\mathbf{\hat\varphi}]_{b_i},~[\mathbf{\tilde\varphi}_{\tau}]_{c_i},
~[\mathbf{\hat\varphi}_{\tau}]_{d_i}),
\texttt{for all }t\in \mathbb{R}\setminus\{0\},\\
c\hat\varphi'_i(t)-\frac{d_i}{c^2}\hat\varphi''_i(t)\leq
f_i(\hat\varphi_i,~[\mathbf{\hat\varphi}]_{a_i},~[\mathbf{\tilde\varphi}]_{b_i},
~[\mathbf{\hat\varphi}_{\tau}]_{c_i},
~[\mathbf{\tilde\varphi}_{\tau}]_{d_i}), \texttt{for all }t\in
\mathbb{R}\setminus\{0\}.
\end{array}\label{b17}
\ees for $i=1,  \cdots, n$. Then $\mathbf{\tilde\varphi}$ and
$\mathbf{\hat\varphi}$ are called coupled quasi-upper solution and
quasi-lower solution of (\ref{b2}), respectively.
\end{defi}

\begin{lem}
If $\mathbf{\tilde\varphi}$, $\mathbf{\hat\varphi}\in
C_{[0,\mathbf{K}]}(\mathbb{R},\mathbb{R}^n)$ are a pair of
quasi-upper and quasi-lower solutions of (\ref{b2}). Then
$\mathbf{\overline x}^{(1)}$ and $\mathbf{\underline x}^{(1)}$
constructed by (\ref{b10}), where $\mathbf{\tilde\varphi}$ and
$\mathbf{\hat\varphi}$ are replaced with these quasi-upper and
quasi-lower solutions, are a pair of upper and lower solutions of
(\ref{b2}), moreover \bess \mathbf{\hat
\varphi}\leq\mathbf{\underline x}^{(1)}\leq
  \mathbf{\overline x}^{(1)}\leq\mathbf{\tilde \varphi}.
\eess\label{lemc2}
\end{lem}
\begin{proof}
Combining (\ref{b17}) and (\ref{b10}) yields \bess\begin{split}
(\overline x_i^{(1)})'-\frac{d_i}{c^2}(\overline
x_i^{(1)})''+\beta_i\overline x_i^{(1)}&=\beta_i\tilde \varphi _i+
f_i(\tilde\varphi_i, [\mathbf{\tilde\varphi}]_{a_i}, [\mathbf{\hat\varphi}]_{b_i}, [\mathbf{\tilde\varphi}_{\tau}]_{c_i}, [\mathbf{\hat\varphi}_{\tau}]_{d_i})\\
&\leq \tilde\varphi'_i-\frac{d_i}{c^2}\tilde\varphi''_i+\beta_i
\tilde\varphi_i,\texttt{for all }t\in \mathbb{R}\setminus\{0\}.
\end{split}\eess
In view of Lemma \ref{lemc1}, we have \bes\begin{split} \overline
x_i^{(1)}(t)&\leq\frac{1}{d_i(\lambda_{2i}-\lambda_{1i})}(\int_{-\infty}^t
e^{\lambda_{1i}(t-s)}(\tilde\varphi'_i(s)-\frac{d_i}{c^2}\tilde\varphi''_i(s)+\beta_i
\tilde\varphi_i(s))ds\\&+\int_{t}^{+\infty}
e^{\lambda_{2i}(t-s)}(\tilde\varphi'_i(s)-\frac{d_i}{c^2}\tilde\varphi''_i(s)+\beta_i
\tilde\varphi_i(s))ds)\\&=\tilde\varphi_i(t),\texttt{ for all }t\in
\mathbb{R}\setminus\{0\}.
\end{split}\label{cc3}\ees
In a similar way, we have $\underline
x_i^{(1)}(t)\geq\hat\varphi_i(t)$ for all $t\in
\mathbb{R}\setminus\{0\}$. Since that $\mathbf{\overline x}^{(1)}$,
$\mathbf{\underline x}^{(1)}$ satisfy (\ref{b10}),
$\mathbf{\overline x}^{(1)}$, $\mathbf{\underline x}^{(1)}\in
\mathcal{C}_{[0, \mathbf{K}]}(\mathbb{R},\mathbb{R}^n)\cap
\mathcal{C}_b^2(\mathbb{R}, \mathbb{R}^n)$ and (\ref{cc3}) holds for
all $t\in\mathbb{R}$. A similar argument in Lemma 2.1 of \cite{lin}
shows that $\mathbf{\overline x}^{(1)}$, $\mathbf{\underline
x}^{(1)}$ are a pair of coupled upper and lower solutions of
(\ref{b2}).
\end{proof}
\begin{theo}
Assume that $(H_1)$ and $(H_2)$ hold. Suppose that
$\mathbf{\tilde\varphi}\in \mathcal{C}_{[0, \mathbf{K}]}(\mathbb{R},
\mathbb{R}^n)$, $\mathbf{\hat\varphi}\in \mathcal{C}_{[0,
\mathbf{K}]}(\mathbb{R}, \mathbb{R}^n)$ be a pair of coupled
quasi-upper and quasi-lower solutions of (\ref{b2}), and \bess
\lim_{t\rightarrow-\infty} \mathbf{\tilde\varphi}(t)=0,\quad
\lim_{t\rightarrow+\infty} \mathbf{\hat\varphi}(t)= \mathbf{K}.\eess
Then (\ref{b2}) and (\ref{b4}) admit a solution. That is, the
problem (\ref{b1}) has a traveling wave solution.
\end{theo}
\begin{proof}
Let \bess\begin{split} \overline
x_i^{(1)}=\frac{c^2}{d_i(\lambda_{2i}-\lambda_{1i})}(\int_{-\infty}^t
e^{\lambda_{1i}(t-s)}(\beta_i\tilde\varphi_i+f_i(\tilde\varphi_i,[\mathbf{\tilde\varphi}]_{a_i},
[\mathbf{\hat\varphi}]_{b_i},
~[\mathbf{\tilde\varphi}_{\tau}]_{c_i},
~[\mathbf{\hat\varphi}_{\tau}]_{d_i}))ds\\+\int_{t}^{+\infty}
e^{\lambda_{2i}(t-s)}(\beta_i\tilde\varphi_i+f_i(\tilde\varphi_i,[\mathbf{\tilde\varphi}]_{a_i},
[\mathbf{\hat\varphi}]_{b_i}, ~[\mathbf{\tilde\varphi}_{\tau}]_{c_i},~[\mathbf{\hat\varphi}_{\tau}]_{d_i}))ds),\\
\underline
x_i^{(1)}=\frac{c^2}{d_i(\lambda_{2i}-\lambda_{1i})}(\int_{-\infty}^t
e^{\lambda_{1i}(t-s)}(\beta_i\hat\varphi_i+f_i(\hat\varphi_i,[\mathbf{\hat\varphi}]_{a_i},
[\mathbf{\tilde\varphi}]_{b_i},
~[\mathbf{\hat\varphi}_{\tau}]_{c_i},~[\mathbf{\tilde\varphi}_{\tau}]_{d_i}))ds\\+\int_{t}^{+\infty}
e^{\lambda_{2i}(t-s)}(\beta_i\hat\varphi_i+f_i(\hat\varphi_i,[\mathbf{\hat\varphi}]_{a_i},
[\mathbf{\tilde\varphi}]_{b_i},
~[\mathbf{\hat\varphi}_{\tau}]_{c_i},~[\mathbf{\tilde\varphi}_{\tau}]_{d_i}))ds)
\end{split}\eess
for $i=1, 2, \cdots, n$. Then by Lemma \ref{lemc2},
$\mathbf{\overline x}^{(1)}$ and $\mathbf{\underline x}^{(1)}$ are a
pair of coupled upper and lower solutions of (\ref{b2}). Replace
$\mathbf{\tilde \varphi}$, $\mathbf{\hat\varphi}$ with
$\mathbf{\overline x}^{(1)}$, $\mathbf{\underline x}^{(1)}$, We then
use Theorem 2.1 and obtain the same results directly.
\end{proof}

\section{Ordered quasi-upper and quasi-lower solutions}
\setcounter{equation}{0}

If $\mathbf{f}(\mathbf{\varphi}, \mathbf{\varphi}_{\mathbf{\tau}})$
is quasimonotone nondecreasing, that is, $b_i=d_i=0$ for all $i$,
then the upper and lower solution are ordered. In this case
\cite{bou,ma,wu} have showed the existence of the traveling wave
solution.

 \begin{defi}
A pair of vectors $\mathbf{\tilde\varphi}\equiv(\tilde
\varphi_1,\cdots,\tilde \varphi_n),$
$\mathbf{\hat\varphi}\equiv(\hat \varphi_1,\cdots,\hat \varphi_n)$
in $\mathcal{C}^2_b(\mathbb{R}, \mathbb{R}^n)$ are called ordered
upper and lower solutions of (\ref{b2}) if
$\mathbf{\tilde\varphi}\geq\mathbf{\tilde\varphi}$ and if \bes
\begin{array}{ll}
\tilde\varphi'_i(t)-\frac{d_i}{c^2}\tilde\varphi''_i(t)\geq
f_i(\mathbf{\tilde\varphi},~\mathbf{\hat\varphi}_{\mathbf{\tau}}),\\
\hat\varphi'_i(t)-\frac{d_i}{c^2}\hat\varphi''_i(t)\leq
f_i(\mathbf{\hat\varphi},~\mathbf{\hat\varphi}_{\mathbf{\tau}}),&(i=1,\cdots,n),
\end{array}\label{e5}
\ees where
$\mathbf{\varphi}_{\mathbf{\tau}}(t)=\mathbf{\varphi}(t-\mathbf{\tau})$.
\end{defi}

\begin{defi}
Assume that $\mathbf{\tilde\varphi}$, $\mathbf{\hat\varphi}\in
\mathcal{C}^1_b(\mathbb{R},\mathbb{R}^n)$ ,
$\frac{d^2\mathbf{\tilde\varphi}}{dt^2}(t)$ and
$\frac{d^2\mathbf{\hat\varphi}}{dt^2}(t)$ exist and continuous on
$\mathbb{R}\setminus\{0\}$, and\bes\begin{split}
\sup_{t\rightarrow\mathbb{R}\setminus\{0\}}
|\frac{d^2\mathbf{\hat\varphi}}{dt^2}(t)|<+\infty, \texttt{ and
}\lim_{t\rightarrow
0^-}\frac{d^2\mathbf{\hat\varphi}}{dt^2}(t),~\lim_{t\rightarrow
0^+}\frac{d^2\mathbf{\hat\varphi}}{dt^2}(t)\texttt{ exist,
}\\
\sup_{t\rightarrow\mathbb{R}\setminus\{0\}}
|\frac{d^2\mathbf{\tilde\varphi}}{dt^2}(t)|<+\infty, \texttt{ and
}\lim_{t\rightarrow
0^-}\frac{d^2\mathbf{\tilde\varphi}}{dt^2}(t),~\lim_{t\rightarrow
0^+}\frac{d^2\mathbf{\tilde\varphi}}{dt^2}(t)\texttt{ exist,
}\label{e6}\end{split}\ees $\mathbf{\tilde\varphi}$,
$\mathbf{\hat\varphi}$ satisfy \bes
\begin{array}{ll}
\tilde\varphi'_i(t)-\frac{d_i}{c^2}\tilde\varphi''_i(t)\geq
f_i(\mathbf{\tilde\varphi},~\mathbf{\tilde\varphi}_{\mathbf{\tau}}),
\texttt{for all }t\in \mathbb{R}\setminus\{0\},\\
\hat\varphi'_i(t)-\frac{d_i}{c^2}\hat\varphi''_i(t)\leq
f_i(\mathbf{\hat\varphi},~\mathbf{\hat\varphi}_{\mathbf{\tau}}),
\texttt{for all }t\in \mathbb{R}\setminus\{0\},
\end{array}\label{e7}
\ees for $i=1,  \cdots, n$. Then $\mathbf{\tilde\varphi}$ and
$\mathbf{\hat\varphi}$ are called ordered quasi-upper solution and
quasi-lower solution of (\ref{b2}), respectively.
\end{defi}

Next theorem shows that if the upper or quasi-upper solution
$\mathbf{\tilde\varphi}(t)$ is nondecreasing respect to $t$, then
the solution $\mathbf{\varphi^*}(t)$ is also nondecreasing. Using
the similar argument in Lemma 4.1 of \cite{lin}, we can induce the
existence of the traveling wavefront of (\ref{b2}).
\begin{theo}
Assume that $(H_1)$ and $(H_2)$ hold and
$\mathbf{f}(\mathbf{\varphi}, \mathbf{\varphi}_{\mathbf{\tau}})$ is
quasimonotone nondecreasing. Suppose that $\mathbf{\tilde\varphi}\in
\mathcal{C}_{[0, \mathbf{K}]}(\mathbb{R}, \mathbb{R}^n)$,
$\mathbf{\hat\varphi}\in \mathcal{C}_{[0, \mathbf{K}]}(\mathbb{R},
\mathbb{R}^n)$ be a pair of ordered quasi-upper and quasi-lower
solutions of (\ref{b2}), and $\mathbf{\tilde\varphi}(t)$ is
nondecreasing with respect to $t$, \bes \lim_{t\rightarrow-\infty}
\mathbf{\tilde\varphi}(t)=0,\quad \lim_{t\rightarrow+\infty}
\mathbf{\hat\varphi}(t)= \mathbf{K}.\label{e8}\ees
 Then the solution of (\ref{b2}) and (\ref{b4}) $\mathbf{\varphi^*}$ is
nondecreasing with respect to $t$. That is, the problem (\ref{b1})
at least has a traveling wavefront solution.\end{theo} In Theorem
4.1, the condition (\ref{e8}) can be replaced by further
restrictions on $\mathbf{f}$, that is \bess (H_1^*)&
~\mathbf{f}(\mathbf{u},\mathbf{u_\tau})|_{\mathbf{u}=0}=\mathbf{f}(\mathbf{u},\mathbf{u_\tau})|_{\mathbf{u}=\mathbf{K}}=0
\texttt{ and
}\mathbf{f}(\mathbf{u},\mathbf{u_\tau})|_{\mathbf{u}=\mathbf{L}}\neq0
\texttt{ for the }\\&  \texttt{constant-valued function } \mathbf{L}
\texttt{ with } 0\leq\mathbf{L}\leq\mathbf{K} \texttt{ and }
 \mathbf{L}\neq0,
\mathbf{L}\neq\mathbf{K}.\eess

As similar as the argument in Theorem 4.4 of \cite{lin}, Theorem 4.1
can be transformed into

\begin{theo}
Assume that $(H_1^*)$ and $(H_2)$ hold and
$\mathbf{f}(\mathbf{\varphi}, \mathbf{\varphi}_{\mathbf{\tau}})$ is
quasimonotone nondecreasing. Suppose that $\mathbf{\tilde\varphi}\in
\mathcal{C}_{[0, \mathbf{K}]}(\mathbb{R}, \mathbb{R}^n)$,
$\mathbf{\hat\varphi}\in \mathcal{C}_{[0, \mathbf{K}]}(\mathbb{R},
\mathbb{R}^n)$ be a pair of ordered quasi-upper and quasi-lower
solutions of (\ref{b2}), and $\mathbf{\tilde\varphi}(t)$ is
nondecreasing with respect to $t$,\bes
0\leq\mathbf{\hat\varphi}(t)\leq\mathbf{\tilde\varphi}(t)\leq
\mathbf{K},\quad
\mathbf{\hat\varphi}(t)\neq0,\mathbf{\tilde\varphi}(t)\neq
\mathbf{K} \texttt{ in } \mathbb{R},\label{e9}\ees
 then the solution of (\ref{b2}) and (\ref{b4}) $\mathbf{\varphi^*}$ is
nondecreasing with respect to $t$. That is, the problem (\ref{b1})
at least has a traveling wavefront solution.\end{theo}

\begin{rem} In Theorems  4.1 and 4.2, the ordered quasi-lower solution is not necessary nondecreasing.
 The results of Theorem 4.1  has been obtained by many authors, for example, Theorem 2.2 of \cite{ma} and Theorem 11 of \cite{bou}, and
 Theorem 3.6 of \cite{wu}.
\end{rem}

In the similar way, Theorem 3.1 can be transformed into

\begin{theo}
Assume that $(H_1^*)$ and $(H_2)$ hold. Suppose that
$\mathbf{\tilde\varphi}\in \mathcal{C}_{[0, \mathbf{K}]}(\mathbb{R},
\mathbb{R}^n)$, $\mathbf{\hat\varphi}\in \mathcal{C}_{[0,
\mathbf{K}]}(\mathbb{R}, \mathbb{R}^n)$ be a pair of coupled
quasi-upper and quasi-lower solutions of (\ref{b2}), and \bes
0\leq\mathbf{\hat\varphi}(t)\leq\mathbf{\tilde\varphi}(t)\leq
\mathbf{K},\quad
\mathbf{\hat\varphi}(t)\neq0,\mathbf{\tilde\varphi}(t)\neq
\mathbf{K} \texttt{ in } \mathbb{R},\label{e10}\ees Then (\ref{b2})
and (\ref{b4}) admit a solution. That is, the problem (\ref{b1}) has
a traveling wave solution.
\end{theo}

\section{Application}
\setcounter{equation}{0}

In this section, we use our results obtained in previous sections to
consider the delayed reaction diffusion models.

\subsection{Belousov-Zhabotinskii equations}
Consider the delayed Belousov-Zhabotinskii equations
 \bes\left\{
\begin{array}{ll}
\frac{\partial u(x, t)}{\partial t}-\frac{\partial^2 u(x,
t)}{\partial x^2}
= u(x, t)(1-u(x, t)-rv(x, t-\tau_2)),&\\
\frac{\partial v(x, t)}{\partial t}-\frac{\partial^2 v(x,
t)}{\partial x^2} =-bu(x,t-\tau_1)v(x,t),&
\end{array}\right.\label{f1}
\ees where $u(x, t), v(x, t)$ are scalar functions, $r>0,b>0$ are
constants. In the biological sense, $u$ and $v$ represent the Bromic
acid and bromide ion concentrations respectively (see more details
in \cite{mur}). Without the delays, the existences of the traveling
wave solutions were considered in \cite{kan,kap,tro,ye}. When
$\tau_1=0, \tau_2\neq0,$ \cite{bou,ma,wu} studied the traveling wave
solution by using of various methods to construct quasi-upper
solution.

Now we seek the traveling wave solution, whose form is
$u(x,t)=\varphi_1(t+x/c)$, $v(x,t)=\varphi_2(t+x/c)$. Clearly the
wave equations corresponding to (\ref{b2}) is the following form
\bes \left\{\begin{array}{l}
              \varphi_1'(t)-\frac{1}{c^2} \varphi_1''(t)=
\varphi_1(t)(1-\varphi_1(t)-r\varphi_2(t-\tau_2)),\\
 \varphi_2'(t)-\frac{1}{c^2} \varphi_2''(t)=-b
\varphi_2(t)\varphi_1(t-\tau_1).\end{array}\right.\label{f2} \ees We
seek a traveling wave solution of (\ref{f1}) with the boundary
conditions \bess \lim_{t\rightarrow-\infty}\varphi_1(t)=0,\quad
\lim_{t\rightarrow-\infty}\varphi_2(t)=0,\\
\lim_{t\rightarrow+\infty}\varphi_1(t)=1,\quad
\lim_{t\rightarrow+\infty}\varphi_2(t)=1.\\
\eess It is easy to check that $(H_1^*)$ and $(H_2)$ are satisfied.
We only need seek the coupled quasi-upper and quasi-lower solutions
of (\ref{f2}).

If $\mathbf{\tilde\varphi}=(\tilde\varphi_1,\tilde\varphi_2)$ and
$\mathbf{\hat\varphi}=(\hat\varphi_1,\hat\varphi_2)$ are coupled
quasi-upper and quasi-lower solutions of (\ref{f2}), they must
satisfy \bes
\begin{array}{ll}
\tilde\varphi'_1(t)-\frac{1}{c^2}\tilde\varphi''_1(t)\geq
\tilde\varphi_1(t)(1-\tilde\varphi_1(t)-r\hat\varphi_2(t-\tau_2)),
\texttt{for all }t\in \mathbb{R}\setminus\{0\},\\
\tilde\varphi'_2(t)-\frac{1}{c^2}\tilde\varphi''_2(t)\geq
-b\tilde\varphi_2(t)\hat\varphi_1(t-\tau_1),
\texttt{for all }t\in \mathbb{R}\setminus\{0\},\\
\hat\varphi'_1(t)-\frac{1}{c^2}\hat\varphi''_1(t)\leq
\hat\varphi_1(t)(1-\hat\varphi_1(t)-r\tilde\varphi_2(t-\tau_2)),
\texttt{for all }t\in
\mathbb{R}\setminus\{0\},\\
\hat\varphi'_2(t)-\frac{1}{c^2}\hat\varphi''_2(t)\leq
-b\hat\varphi_2(t)\tilde\varphi_1(t-\tau_1),
\texttt{for all }t\in \mathbb{R}\setminus\{0\}.\\
\end{array}\label{f3}
\ees Assume that $\mathbf{\tilde\varphi}$ and
$\mathbf{\hat\varphi}$\bes
\tilde\varphi_1(t)=\left\{\begin{array}{ll}
       \frac{1}{2}e^{\lambda_1t},&t\leq0,\\
 1-\frac{1}{2}e^{-\lambda_1t},&t>0,\end{array}\right. \texttt{ and } \tilde\varphi_2(t)=\left\{\begin{array}{ll}
       \frac{1}{2}e^{\lambda_2t},&t\leq0,\\
 1-\frac{1}{2}e^{-\lambda_2t},&t>0,\end{array}\right.\label{f4}\ees
 \bes
\hat\varphi_1(t)=\left\{\begin{array}{ll}
       \delta ke^{\lambda_3t},&t\leq0,\\
 k-\delta k e^{-\lambda_3 t},&t>0,\end{array}\right. \texttt{ and } \hat\varphi_2(t)=0,\label{f5}\ees
 where $\lambda_1,\lambda_2,\lambda_3,\delta,k$ are undetermined positive
 constants.

 Direct calculations show that
 \bes
\tilde\varphi'_1(t)=\left\{\begin{array}{ll}
       \frac{\lambda_1}{2}e^{\lambda_1t},&t\leq0,\\
 \frac{\lambda_1}{2}e^{-\lambda_1t},&t>0,\end{array}\right. \texttt{ and } \tilde\varphi'_2(t)=\left\{\begin{array}{ll}
       \frac{\lambda_2}{2}e^{\lambda_2t},&t\leq0,\\
 \frac{\lambda_2}{2}e^{-\lambda_2t},&t>0,\end{array}\right.\label{f6}\ees
 \bes
\tilde\varphi''_1(t)=\left\{\begin{array}{ll}
       \frac{\lambda_1^2}{2}e^{\lambda_1t},&t\leq0,\\
- \frac{\lambda_1^2}{2}e^{-\lambda_1t},&t>0,\end{array}\right.
\texttt{ and } \tilde\varphi''_2(t)=\left\{\begin{array}{ll}
       \frac{\lambda_2^2}{2}e^{\lambda_2t},&t\leq0,\\
 -\frac{\lambda_2^2}{2}e^{-\lambda_2t},&t>0,\end{array}\right.\label{f7}\ees
 \bes
\hat\varphi'_1(t)=\left\{\begin{array}{ll}
      \lambda_3 \delta ke^{\lambda_3t},&t\leq0,\\
\lambda_3\delta k e^{-\lambda_3 t},&t>0,\end{array}\right. \texttt{
and } \hat\varphi''_1(t)=\left\{\begin{array}{ll}
      \lambda_3^2 \delta ke^{\lambda_3t},&t\leq0,\\
-\lambda_3^2\delta k e^{-\lambda_3
t},&t>0,\end{array}\right.\label{f8}\ees From (\ref{f6}),(\ref{f7})
and (\ref{f8}), we see that the first derivatives are continuous and
the second derivatives exist and continuous on
$\mathbb{R}\setminus\{0\}$. Hence $\mathbf{\tilde\varphi}$ and
$\mathbf{\hat\varphi}$ satisfy (\ref{e6}). Now we will choose proper
$\lambda_1,\lambda_2,\lambda_3,\delta,k$
 such that (\ref{f3}) holds.

Substituting (\ref{f4}), (\ref{f6}) and (\ref{f7}) into the first
equation of (\ref{f3}) yields \bes\begin{split}
&\frac{\lambda_1}{2}e^{\lambda_1t}-\frac{1}{c^2}\frac{\lambda_1^2}{2}e^{\lambda_1t}\geq
\frac{\lambda_1}{2}e^{\lambda_1t}(1-\frac{1}{2}e^{\lambda_1t})
\texttt{ for }t<0,\\&
\frac{\lambda_1}{2}e^{-\lambda_1t}-\frac{1}{c^2}(-\frac{\lambda_1^2}{2}e^{-\lambda_1t})\geq
(1-\frac{1}{2}e^{-\lambda_1 t})(1-1+\frac{1}{2}e^{-\lambda_1 t})
\texttt{ for }t\geq0.
\end{split}\label{f9} \ees In order to induce
(\ref{f9}), the following is sufficient \bes
\lambda_1-\frac{\lambda_1^2}{c^2}-1\geq0.\label{f10} \ees Hence we
set \bes
\lambda_1=\frac{c^2(1-\sqrt{1-\frac{4}{c^2}})}{2}.\label{f11} \ees

Substituting (\ref{f4})-(\ref{f8}) into the second equation of
(\ref{f3}) yields \bes\begin{split}
&\frac{\lambda_2}{2}e^{\lambda_2t}-\frac{1}{c^2}\frac{\lambda_2^2}{2}e^{\lambda_2t}\geq
-b\frac{\lambda_1}{2}e^{\lambda_2t}\delta k e^{\lambda_3(t-\tau_1)}
\texttt{ for
}t<0,\\
&\frac{\lambda_2}{2}e^{-\lambda_2t}-\frac{1}{c^2}(-\frac{\lambda_2^2}{2}e^{-\lambda_2t})\geq
-b(1-\frac{1}{2}e^{-\lambda_2t})\delta k e^{\lambda_3(t-\tau_1)}
\texttt{ for }0\leq t\leq \tau_1,
\\&
\frac{\lambda_2}{2}e^{-\lambda_2t}-\frac{1}{c^2}(-\frac{\lambda_2^2}{2}e^{-\lambda_2t})\geq
-b(1-\frac{1}{2}e^{-\lambda_2 t})(k-\delta
ke^{-\lambda_3(t-\tau_1)}) \texttt{ for }t>\tau_1.
\end{split}\label{f12} \ees In order to induce
(\ref{f12}), the following is sufficient \bes
\lambda_2-\frac{\lambda_2^2}{c^2}\geq0.\label{f13} \ees Hence we set
\bes \lambda_2=\varepsilon_1, \texttt{ where
}\varepsilon_1<<1.\label{f14} \ees

Substituting (\ref{f4})-(\ref{f8}) into the third equation of
(\ref{f3}) yields \bes\begin{split} &\lambda_3\delta ke^{\lambda_3
t}-\frac{1}{c^2}\lambda_3^2\delta k e^{\lambda_3 t}\leq \delta k
e^{\lambda_3t}(1-\delta ke^{\lambda_3
t}-r\frac{1}{2}e^{\lambda_2(t-\tau_2)}) \texttt{ for
}t<0,\\
&\lambda_3\delta ke^{-\lambda_3 t}+\frac{1}{c^2}\lambda_3^2\delta k
e^{-\lambda_3 t}\leq (k-\delta ke^{-\lambda_3 t})(1-k+\delta
ke^{-\lambda_3 t}-r\frac{1}{2}e^{\lambda_2(t-\tau_2)}) \\&
\texttt{
for }0\leq t\leq \tau_2,
\\&
\delta k e^{-\lambda_3t}(\lambda_3+\frac{\lambda_3^2}{c^2})\leq
(k-\delta ke^{-\lambda_3 t})(1-k+\delta ke^{-\lambda_3
t}-r(1-\frac{1}{2}e^{-\lambda_2t}))
\\&\texttt{ for }t>\tau_2.
\end{split}\label{f15} \ees
If we set $\delta<<1,\lambda_2<<1$, in order to induce (\ref{f15}),
the following is sufficient \bes\begin{split}&
\lambda_3-\frac{\lambda_3^2}{c^2}-1<0,&\\&\delta k e^{-\lambda_3
t}(\lambda_3+\frac{\lambda_3^2}{c^2}+1-k)<k(1-k),&
\\&
 \delta
ke^{-\lambda_3
t}(\lambda_3+\frac{\lambda_3^2}{c^2}+1-r-k)<k(1-r-k).&
\end{split}\label{f16}\ees  To ensure (\ref{f16}), the
parameter must satisfy \bes r<1 \label{f17}\ees and we set \bes
\lambda_3=\frac{c^2(1-\sqrt{1-\frac{4}{c^2}})}{2}-\varepsilon_2,
\texttt{ where } \varepsilon_2<<1,k<<1.\ees It is easy to see that
the fourth equation of (\ref{f3}) naturally holds. Therefore, we
have proved
\begin{lem} Assume that the parameter of the problem (\ref{f1})
$r<1$. Then there exists a constant $c^*=2$ such that if $c>c^*$,
$\mathbf{\tilde\varphi}$ and $\mathbf{\hat\varphi}$, which are
defined in (\ref{f4}) and  (\ref{f4}), are a pair of coupled
quasi-upper and quasi-lower solutions of the problem
(\ref{f2}).\label{lemd1}
\end{lem}

Finally, by Theorem 4.3 we have
\begin{theo} Assume that the parameter of the problem (\ref{f1})
$r<1$. Then there exists a constant $c^*=2$ such that if $c>c^*$,
the problem (\ref{f1}) has a traveling wave solution
$u(x,t)=\varphi_1(t+x/c)$, $v(x,t)=\varphi_2(t+x/c)$, which connects
$(0,0)$ and $(1,1)$.\label{lemd2}
\end{theo}

\begin{rem} In \cite{bou}, the critical value of wave velocity $c^*$
is dependent on the parameter $b$. In our Theorem 5.1, we show that
$c^*$ may be a constant.
\end{rem}

\subsection{Mutualistic Lotka-Volterra model}
The delayed mutualistic Lotka-Volterramodel is as follows
\bes\left\{
\begin{array}{ll}
\frac{\partial u(x, t)}{\partial t}-d_1\frac{\partial^2 u(x,
t)}{\partial x^2}
= r u(x, t)(1-a_1u(x,t)+b_1v(x, t)),&\\
\frac{\partial v(x, t)}{\partial t}-d_2\frac{\partial^2 v(x,
t)}{\partial x^2} =v(x,t)(a_2u(x,t-\tau)-b_2),&
\end{array}\right.\label{g1}
\ees where $u(x,t)$, $v(x,t)$ are scalar functions, and
$r,a_1,a_2,b_1,b_2$ are all positive constants, $d_1,d_2$ are the
positive diffusion coefficients. $\tau$ represent the positive
delay. For a detailed description of this model, we refer to
\cite{mur}. The above model with $\tau=0, b_1<0$ has been considered
in \cite{dun1,dun2,mur}. The case with $\tau=0, b_1<0$ was studied
in \cite{bou,ma}. However, when $b_1<0$ the reaction term of
(\ref{g1}) is not quasimonotone nondecreasing,  there was not
sufficient condition to prove Lemma 3.5 of \cite{ma} and Lemma 12 of
\cite{bou}.

If we assume that \bes a_2>a_1b_2, \label{g2}\ees then the model
(\ref{g1}) has a unique positive equilibrium \bes
(u^*,v^*)=(\frac{b_2}{a_2},\frac{1}{b_1}(\frac{a_1b_2}{a_2}-1)).\label{g3}
\ees

The wave equations corresponding to (\ref{b2}) is the following form
\bes \left\{\begin{array}{l}
              \varphi_1'(t)-\frac{d_1}{c^2} \varphi_1''(t)=r
\varphi_1(t)(1-a_1\varphi_1(t)+b_1\varphi_2(t)),\\
 \varphi_2'(t)-\frac{d_2}{c^2} \varphi_2''(t)=
\varphi_2(t)(a_2\varphi_1(t-\tau)-b_2).\end{array}\right.\label{g4}
\ees We will find the solution of the above wave equation such that
\bess \lim_{t\rightarrow-\infty}\varphi_1(t)=0,\quad
\lim_{t\rightarrow-\infty}\varphi_2(t)=0,\\
\lim_{t\rightarrow+\infty}\varphi_1(t)=u^*,\quad
\lim_{t\rightarrow+\infty}\varphi_2(t)=v^*.\\
\eess

It is also easy to verify that $(H_1^*)$ and $(H_2)$ are satisfied.
We only need seek the coupled quasi-upper and quasi-lower solutions
of (\ref{g4}).

If $\mathbf{\tilde\varphi}=(\tilde\varphi_1,\tilde\varphi_2)$ and
$\mathbf{\hat\varphi}=(\hat\varphi_1,\hat\varphi_2)$ are coupled
quasi-upper and quasi-lower solutions of (\ref{g4}), they must
satisfy \bes
\begin{array}{ll}
\tilde\varphi'_1(t)-\frac{d_1}{c^2}\tilde\varphi''_1(t)\geq
r\tilde\varphi_1(t)(1-a_1\tilde\varphi_1(t)+b_1\tilde\varphi_2(t)),
\texttt{for all }t\in \mathbb{R}\setminus\{0\},\\
\tilde\varphi'_2(t)-\frac{d_2}{c^2}\tilde\varphi''_2(t)\geq
\tilde\varphi_2(t)(a_2\tilde\varphi_1(t-\tau)-b_2),
\texttt{for all }t\in \mathbb{R}\setminus\{0\},\\
\hat\varphi'_1(t)-\frac{d_1}{c^2}\hat\varphi''_1(t)\leq
r\hat\varphi_1(t)(1-a_1\hat\varphi_1(t)-b_1\hat\varphi_2(t)),
\texttt{for all }t\in
\mathbb{R}\setminus\{0\},\\
\hat\varphi'_2(t)-\frac{d_2}{c^2}\hat\varphi''_2(t)\leq
\hat\varphi_2(t)(a_2\hat\varphi_1(t-\tau)-b_2),
\texttt{for all }t\in \mathbb{R}\setminus\{0\}.\\
\end{array}\label{g5}
\ees

Assume that $\mathbf{\tilde\varphi}$ and $\mathbf{\hat\varphi}$\bes
\tilde\varphi_1(t)=\left\{\begin{array}{ll}
       \frac{u^*}{2}e^{\lambda_1t},&t\leq0,\\
u^*-\frac{u^*}{2}e^{-\lambda_1t},&t>0,\end{array}\right. \texttt{
and } \tilde\varphi_2(t)=\left\{\begin{array}{ll}
       \frac{v^*}{2}e^{\lambda_2t},&t\leq0,\\
v^*-\frac{v^*}{2}e^{-\lambda_2t},&t>0,\end{array}\right.\label{g6}\ees
 \bes
\hat\varphi_1(t)=\left\{\begin{array}{ll}
       \delta ke^{\lambda_3t},&t\leq0,\\
 k-\delta k e^{-\lambda_3 t},&t>0,\end{array}\right. \texttt{ and } \hat\varphi_2(t)=0,\label{g7}\ees
 where $\lambda_1,\lambda_2,\lambda_3,\delta,k$ are undetermined positive
 constants.

 As similar as the process in the above subsection, in order to
 (\ref{g5}), the following is sufficient

 \bes\begin{split}
&\lambda_1-\frac{d_1\lambda_1^2}{c^2}-r\geq0,\\&
\lambda_1+\frac{d_1\lambda_1^2}{c^2}-b_1v^*\geq0,\\
& \lambda_2(1-\frac{d_2}{c^2}\lambda_2)\geq0,\\
& \lambda_3-\frac{\lambda_3^2}{c^2}\leq r,\\
& k<<1,\delta<<1.
\end{split}\label{g8} \ees

If \bes c>\max\{2\sqrt{rd_1},\sqrt{d_1(r+\frac{a_1b_2}{a_2}-1)}\},
\label{g9} \ees set\bes
\lambda_1=\frac{c^2(1+\sqrt{1-\frac{4rd-1}{c^2}})}{2d_1},\lambda_2<<1,\lambda_3<<1,k<<1,\delta<<1,
\label{g10} \ees then the sufficient conditions (\ref{g8}) hold.
Therefore by Theorem 4.3, we have
\begin{theo} Assume that the parameter of the problem (\ref{g1})
satisfy (\ref{g2}). Then there exists a constant
$c^*=\max\{2\sqrt{rd_1},\sqrt{d_1(r+\frac{a_1b_2}{a_2}-1)}\}$ such
that if $c>c^*$, the problem (\ref{f1}) has a traveling wave
solution $u(x,t)=\varphi_1(t+x/c)$, $v(x,t)=\varphi_2(t+x/c)$, which
connects $(0,0)$ and
$(\frac{b_2}{a_2},\frac{1}{b_1}(\frac{a_1b_2}{a_2}-1))$.\label{lemd3}
\end{theo}

\section{Discussion}

We aim to study the existence for the traveling wave solution of the
discrete-delayed reaction diffusion systems, where the reaction term
is mixed quasimonotone. Our result is that the existence of the
coupled quasi-upper and quasi-lower solutions ensure the traveling
wave solution exists. In the equations of \cite{bou,ma,wu}, the
reaction term is quasimonotone nondecreasing. In fact the conditions
of quasimonotone nondecreasing property is very strong. The
predator-prey model which is studied in \cite{bou,ma} do not satisfy
the quasimonotone nondecreasing property, and the existence theorem
of \cite{bou,ma} is not suitable. Hence the application scope of
Theorem 3.1 is wider than \cite{bou,ma,wu}. Moreover, as a special
case of Theorem 3.1, Theorem 4.1 contains the previous existence
theorem for the traveling wave solution.

Our technique to deal with the mixed quasimonotonity is constructing
the coupled upper and lower solutions. Recently the method so-called
cross iteration scheme was developed in \cite{li} to deal with the
traveling wave solution for the 2 dimensional competitive
Lotka-Volterra model. Comparing with the model of \cite{li}, the
systems in this paper are extensive to $n$ dimension.

The classical coupled upper and lower solutions need the second
order smoothness. It is very difficult to satisfy this condition in
the real model. Thus it is necessary to relax the smoothness of the
coupled upper and lower solutions to first order smoothness, which
is called coupled quasi-upper and quasi-lower solutions. This paper
apply the modified Perron theorem with the case $C^1$ smoothness,
which is first proposed in \cite{bou}. Our existence theorem of the
traveling wave solution is suitable to all 2 species or 3 species
Lotka-Volterra systems.

\end{document}